\documentclass[aps,rmp]{revtex4}
\usepackage{amsfonts,amssymb,amsmath,bm}
\usepackage{graphics,epsfig}
\usepackage{caption}
\usepackage{verbatim,ulem}
\usepackage{hyperref}
\usepackage{breakurl}
\usepackage{here}
\usepackage{color}
\usepackage[ngerman,french,english]{babel}

\newcommand{\ue}{\mathrm{e}}

\def\be{\begin{eqnarray}}
\def\en{\end{eqnarray}}
\def\UB{\begin{equation}}
\def\UE#1{\label{#1}\end{equation}}
\def\BE{\begin{equation}}
\def\EE#1{\label{#1}\end{equation}}

\def\R{{\mathbb R}}

\def\nL{\nabla^{\rm L}}
\def\rf#1{(\ref{#1})}

\def\bv{{\bm v}}
\def\bx{{\bm x}}

\def\ba{{\bm a}}
\def\bF{{\bm F}}
\def\bn{{\bm n}}
\def\bom{{\bm\omega}}

\begin{document}
\title{Cauchy's almost forgotten Lagrangian
formulation of the Euler equation for 3D incompressible flow}
\author{Uriel Frisch}
\affiliation{UNS,~CNRS,~OCA,~Lab.~Lagrange,~B.P.~4229,~06304~Nice~Cedex~4,~France}
\author{Barbara Villone}
\affiliation{INAF, Osservatorio Astrofisico di Torino, Via Osservatorio, 20, 10025
Pino Torinese, Italy}
\date{\today}
\begin{abstract}
Two prized papers, one by Augustin Cauchy in 1815, presented to the
French Academy and the other by Hermann Hankel in 1861, presented
to G\"ottingen University, contain major discoveries on vorticity
dynamics whose impact is now quickly increasing. Cauchy found a
Lagrangian formulation of 3D ideal incompressible flow in terms of
three invariants that generalize to three dimensions the now well-known law of 
conservation of vorticity along fluid particle trajectories for
two-dimensional flow.  This has very recently been used to prove analyticity in time of fluid
particle trajectories for 3D incompressible Euler flow and can be
extended to compressible flow, in particular to cosmological dark matter. Hankel showed that Cauchy's
formulation gives a very simple Lagrangian derivation of the Helmholtz
vorticity-flux invariants and, in the middle of the proof, derived an
intermediate result which is the conservation of the circulation of
the velocity around a closed contour moving with the fluid. This
circulation theorem was to be rediscovered independently by William
Thomson (Kelvin) in 1869.  Cauchy's invariants were only occasionally
cited in the 19th century --- besides Hankel, foremost by George Stokes and Maurice
L\'evy --- and even less so in the 20th until they were rediscovered via
Emmy Noether's theorem in the late 1960, but reattributed to Cauchy only
at the end of the 20th century by Russian scientists.
\end{abstract}
\maketitle

\section{Introduction}
\label{intro}

The motion of a fluid can be described either in a fixed frame of
reference in terms of the spatio-temporal \textit{Eulerian}
coordinates
$x,\,y,\,z,\,t$ or by following individual fluid particles in terms of their
initial or \textit{Lagrangian} coordinates $a,\,b,\,c,\,t$. The
description of a viscous fluid is simpler in Eulerian coordinates, but 
for an ideal (inviscid) incompressible fluid there is a variational
formulation in Lagrangian coordinates of the equations of motion, due
precisely to Lagrange. Over most of the 19th and 20th century the
Eulerian formulation, which is better adapted to cases where the
boundaries of the flow are prescribed and fixed, became
predominant. Nevertheless, the Lagrangian formulation is more natural
for wave motion in a free-surface flow and for complete disruptions of
the flow, such as the breakup of a dam.\footnote{Lagrange, 1788.
  For more context on the early 18th century developments
in fluid dynamics, see Truesdell, 1954b;
% the Trusdell 1954 reference of Darrigo-Frisch
 Darrigol, 2005;
% the Darrigol 2005 book cited in Darrigol-Frisch
 Darrigol \& Frisch, 2008. 
Regarding the history of Eulerian and Lagrangian coordinates,
  see Truesdell, 1954a. For wave motion, see \S\,\ref{ss:prize}; for
  dam  breakup,
  see the end of \S\,\ref{ss:noether}.
% The 1954 book 
} 

In more recent years there has been a strong renewal of interest in
Lagrangian approaches. In \textit{cosmology} most of the matter is
generally believed to be of the dark type, which is essentially
collisionless and thus inviscid. Lagrangian perturbations methods have
been developed since the nineties that shed light on the mechanisms of
formation of cosmological large-scale structures. Closely related is
the problem of reconstruction of the past Lagrangian history of the
Universe from present-epoch observations of the distribution of
galaxies.
Novel fast
photography techniques have been developed for the tracking of many
particles seeded into laboratory flow that allow the reconstruction of a
substantial fraction of the Lagrangian map that associates the
positions of fluid particles at some initial time to their positions
at later times. When the flow is electrically highly conducting and
can support a magnetic field by the magnetohydrodynamic (MHD) dynamo
effect, it has been shown that the long-time fate of such a magnetic
field is connected to the issue of Lagrangian chaos, namely how fast
neighbouring fluid particle trajectories separate in
time.\footnote{For Lagrangian perturbations in cosmology, see, e.g., Moutarde {\it et al.}, 1991; 
Buchert, 1992. For reconstruction of the Universe, see Brenier
\textit{et al.}, 2003 and references therein. For Lagrangian experimental and numerical tracking techniques, see, e.g., Toschi \&
Bodenschatz, 2009. %paper sent to BV on Jan. 5
For MHD, see Arnold, Zeldovich, Ruzmaikin \& Sokoloff, 1981; Vishik, 1989;
%bv fatto ì. UF: dans mon livre
Childress \& Gilbert, 1995. %dans mon livre
}

At a more fundamental level, over 250 years after Euler wrote the
equations governing incompressible ideal three-dimensional (3D) fluid flow
(generally known as the 3D Euler equations), we still do not know if the
solutions remain well-behaved for all times or become singular and dissipate
energy after some finite time; and this even when the initial data are very
smooth, say, analytic. Many attempts have been made to tackle this
problem by numerical simulations in an Eulerian-coordinates
framework, but the problem remains moot. Lagrange himself frequently
preferred Eulerian 
coordinates, although
the variational formulation he gave was formulated in
Lagrangian coordinates: in modern mathematical language, the solutions to the equations
of incompressible  fluid flow are geodesics on the
infinite-dimensional manifold SDiff of volume-preserving Lagrangian maps.
In a Lagrangian framework one is ``riding the eddy'' and does not feel
too much of the possible spatial roughness. Actually, fluid particle trajectories
can be analytic in time even when the flow has only limited spatial
smoothness.\footnote{Euler, 1757. 
%UF: put the original reference and the translation in the EE250
%Proceedings: English translation by U. Frisch ``General principles of the 
%motion of fluids'' 2008 {\it Physica D \bf 237},
%1825--1839. arXiv 0802.2383 [nlin.CD]
On the issue of well-posedeness, see Eyink, Frisch, Moreau and Sobolevsky, 2008.
%G.~Eyink, U.~Frisch, R.~Moreau and A.~Sobolevsky eds. 2008, 
%Euler Equations: 250 Years On, Proceedings of Conference held in Aussois, June
%18--23, 2007, 450 pages {\it Physica D \bf 237}, issues 14--17, 2008. 
 On analyticity of fluid particle trajectories, see Serfati, 1995;
%J. Math. Pures Appl. See references in Zhelig.-Frisch
Shnirelman, 2012;%Global and Stochastic Analysis See references in Zhelig.-Frisch
 ~see also the next paragraph.}

Recently, borrowing some ideas of the cosmological Lagrangian
perturbation theory, one of us and Vladislav Zheligovsky obtained a
form of the incompressible 3D Euler equations in Lagrangian
coordinates, from which simple recursion relations among the temporal
Taylor coefficients of the Lagrangian map can be derived and
analyticity in time of the Lagrangian trajectories can be proved in a
rather elementary way. This Lagrangian approach can also be used to save
significant computer time in high-resolution simulations of 3D ideal
Euler flow. The Lagrangian equations used for this did not seem at
first glance to be
widely known, but after some time spent searching the past scientific
literature, we found that the equations had been derived in 1815 by
Augustin Cauchy in a long memoir that won a prize from the French
Academy.\footnote{On analyticity of the Lagrangian map, see Frisch \&
  Zheligovsky, 2014; Zheligovsky \& Frisch, 2014.  On a novel Lagrangian
  numerical method, see project NELmezzo,
  available on request. Cauchy, 1815/1827.}

At first sight, Cauchy's equations, to be presented in
Section~\ref{ss:15-17} ---  here called `\textit{Cauchy's
invariants equations} to distinguish them from an important 
corollary, \textit{Cauchy's vorticity
formula},  --- seemed hardly cited at all.  %Had Cauchy's result been
%forgotten for nearly two hundred years? 
We then engaged in a much more
systematic search.  The surprising result can be summarized as
follows: in the 19th century, Cauchy's invariants equations are cited
only in a small number of papers, the most important one being by
Hermann Hankel; in the 20th century, Cauchy's result seems almost
completely uncited, except at the very end of the century.

The outline of the paper is as follows.  Section~\ref{s:cauchy} is devoted to
Augustin Cauchy: in Section~\ref{ss:microbiocauchy} we recall a few
biographical elements, in particular those connected to the present
study; Section~\ref{ss:prize} is devoted to the
1815 prize-winning {\it M\'emoire sur la propagation des ondes}; in
Section~\ref{ss:15-17} we analyze in detail the very beginning of its second
part, which contains Cauchy's Lagrangian formulation of the 3D ideal
incompressible equations in terms of what would now be called
invariants, a terminology we adopt here. 
Section~\ref{s:19th} is about 19th century scientists who realized the
importance of Cauchy's Lagrangian formulation. Outstanding, here, is Hermann Hankel
(Section~\ref{ss:hankel}), a German mathematician whose keen interest
in the history
of mathematics allowed him not to miss Cauchy's 1815 work and to
understand in 1861 its
potential in deriving Helmholtz's results on vorticity in a Lagrangian
framework, and discovering on the way  what is known as  Kelvin's
circulation theorem. 
Then, in Section~\ref{ss:stokesetal}, we
turn to the other 19th century scientist who discuss Cauchy's
invariants equations: Foremost George Stokes, then Maurice L\'evy,  Horace Lamb,
Jules Andrade and Paul Appell and to a
few others who mention the equations but may not be aware that they were
obtained by Cauchy: Gustav Kirchhoff and Henri 
Poincar\'e.\footnote{Hankel, 1861. 
  Helmholtz, 1858. Thomson (Lord Kelvin), 1869.}

Then, in Section~\ref{s:20th}, we turn to the 20th century and
beyond. The first part (Section~\ref{ss:noncauchy}) has Cauchy's
invariants equations apparently fallen into  oblivion. In the second part
(Section~\ref{ss:noether} we shall see that, as a consequence of Emmy Noether's
theorem connecting continuous symmetries and invariants, a number of
scientists were able to rederive Cauchy's invariants equations for the
3D ideal incompressible flow, but without being aware of Cauchy's
work. In Section~\ref{ss:rebirth} we shall find that Russian scientists,
followed by  others, reminded us that all this had been started a long
time ago by Cauchy. In Section~\ref{s:conclusion} we make some concluding remarks.\footnote{Noether, 
1918. For the rediscovery of Cauchy, see Abrashkin, Zen'kovich and Yakubovich, 1996; Zakharov \& Kuznetsov, 1997.}

Since one of our key goals is to understand how important work such as
Cauchy's 1815 formulation of the hydrodynamical equations in
Lagrangian coordinates  managed to get nearly lost,
we are obliged to pay
attention to who cites whom. This is a delicate matter, given that
present-day ethical rules of citing definitely did not apply in past
centuries. But without fast communications, a
peer-review system to point out missing references and a much larger
scientific population, the rules had to be
different. Here, we shall do our best to  mention each instance of citation
of previous work by the author being discussed, when such work is
relevant to our paper.

\section{Cauchy}
\label{s:cauchy}

\subsection{Biographical elements}
\label{ss:microbiocauchy}

\begin{figure}
\centerline{\includegraphics[height=4in]{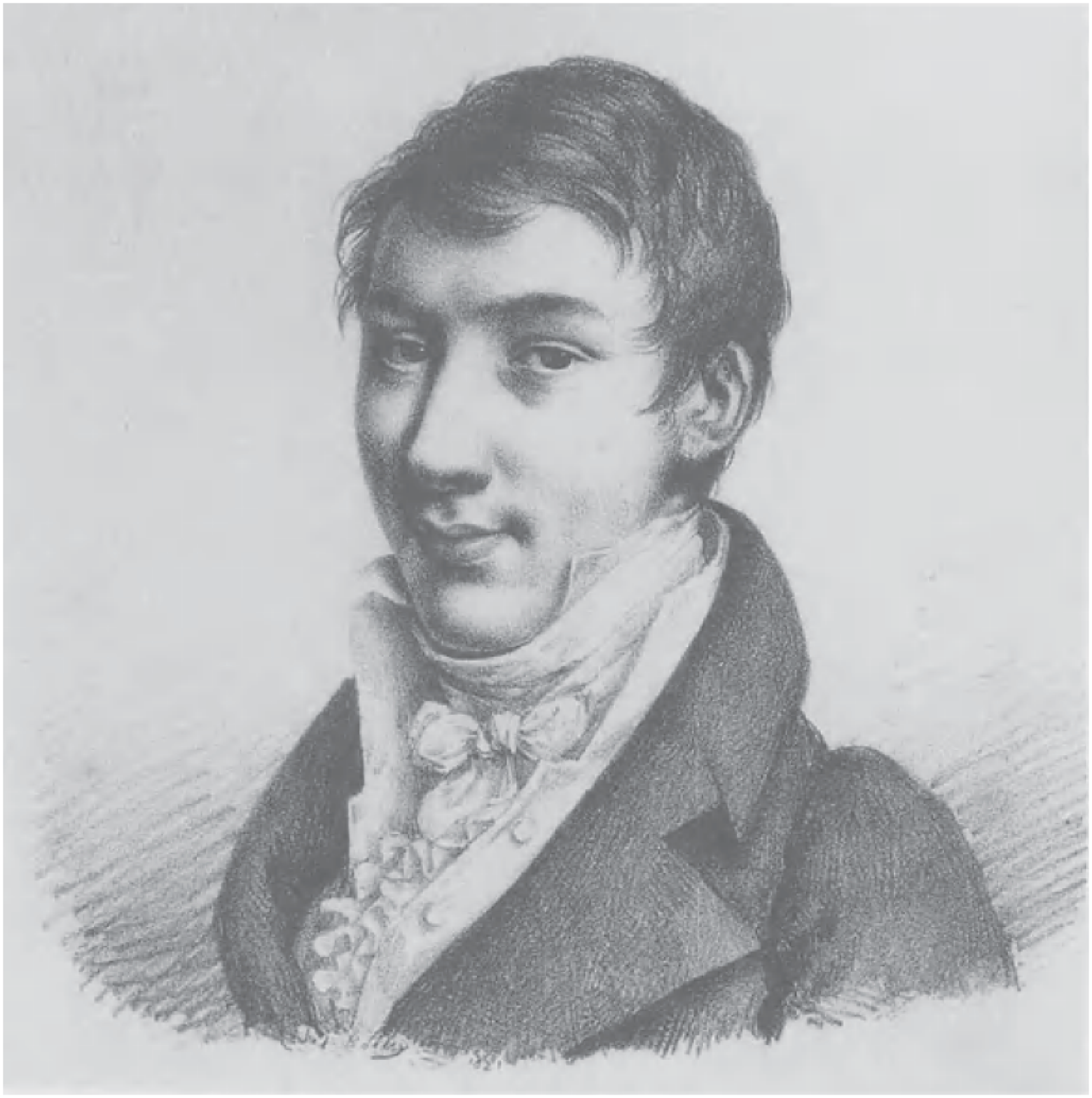}}
\caption{Portrait of Augustin Cauchy by Louis-L\'eopold Boilly (by
  permission 
of \'Ecole Polytechnique).}
\label{fig1}\end{figure}

These  elements are just intended to give a background
on the circumstances which led to Cauchy's 1815 work.  Our main sources have
been the biography of Augustin-Louis Cauchy by Bruno
Belhoste, the biography by Claude Alphonse Valson, written just a few years
after the death of Cauchy and the minutes (\textit{Proc\`es-Verbaux}) of the
meetings of the French Academy of Sciences, referred to as PV, followed by the
date of the corresponding meeting.\footnote{Belhoste, 1991;  cf also his PhD 
thesis, Belhoste, 1982; Valson, 1868; Proc\`es-Verbaux des S\'eances de  l'Acad\'emie, 1915.}

Cauchy was born in 1789, turbulent times, but this did not affect
his ability to get the best theoretical and practical training available in
the early 19th century, attending successively \'Ecole Polytechnique and
\'Ecole des Ponts et Chauss\'ees. His first employment was as a junior
engineer in a major harbour project in Cherbourg in 1810 and then as
an engineer at the Ourcq Canal project in Paris in 1813. During his
engineering years he already displayed keen interest in deep mathematical
questions, several of which he solved in a way that sufficiently impressed the
mathematicians at the Academy of Sciences, (called ``First Class of the
Institute of France'' until King Louis XVIII restored the old name of
``Academy of sciences''; here, we refer to it just as ``Academy''). In
1813, the Geometry section of the Academy ranked Cauchy second for an election to the Academy but the vote
of the rest of the members of the Academy went heavily against him. Anyway,
during the years 1812--1815 Cauchy had a strong coupling to the Academy and his
name appears about one hundred times in the minutes of the Academy
meetings. His awarding of the 1815 mathematics prize by the Academy
(see Section~\ref{ss:prize})
is one more evidence that he was a rising star.\footnote{PV: 1812--1815. For the 1813 election, see  Belhoste, 1991: 37.}

Eventually, the King took advantage of the reorganization of the
Academy to remove Lazare Carnot and Gaspard Monge from the Academy and
appoint Cauchy as a member in March 1816. This unusual way of entering
the Academy produced some friction with regular members.\footnote{Belhoste,
1991: 46--47.}

Cauchy, who was soon to become a world-dominant figure in mathematics
and mathematical physics, was no easy-going personality, he was however 
willing to suffer considerably for his ideas, particularly those
grounded in his Christian beliefs. For example, he went into exile in
1830 for eight long years in order not to have to swear an oath of
allegiance to King Louis Philippe, considered by Cauchy as not
legitimate.\footnote{Belhoste, 1991:  Chaps. 9 and 10.}

\subsection{The 1815 mathematics prize}
\label{ss:prize}

\textit{The 
differential equations given by the author are rigorously applicable only 
to the case where the depth of the fluid is infinite; but he succeeded in 
obtaining their general integrals in a form allowing to discuss the results
and comparing them to experiments.}\footnote{Les \'equations diff\'erentielles
  donn\'ees par l'auteur ne s'appliquent rigoureusement qu'au cas o\`u la
profondeur du fluide est infinie ; mais il est parvenu \`a obtenir leurs
int\'egrales g\'en\'erales sous une forme qui permet d'en discuter les
r\'esultats et de les comparer \`a l'exp\'erience. Institut de France, 1816.}\\[-1.5ex]

Thus reads the beginning of the statement made by the Academy during the public ceremony
of January 8, 1816 at which Cauchy received the 1815 \textit{Grand Prix} 
(mathematics prize). The events leading to this prize started on December 27,
1813 when the Academy committee in charge of proposing a
subject for a mathematics prize decided for ``\textit{le probl\`eme des ondes 
\`a la surface d'un liquide de profondeur ind\'efinie}'' (the problem of waves
on the surface of a liquid of arbitrary depth). The committee put  
Laplace in charge of defining the scientific programme.\footnote{PV: 1815, 27
  d\'ecembre.} 

On October 2, 1815 the Academy received two anonymous manuscripts -- as usual in those circumstances,
distinguished by an epigraph for later identification. Cauchy's manuscript  
had the epigraph  \textit{Nosse quot ionii veniant ad
littora fluctus} (Virgil, Geor.  II, 108; translation: To know how many  
waves come rolling shoreward from the Ionian sea).

On December 26, 1815 the committee proposed giving the 1815 prize to Cauchy's
manuscript.\footnote{PV: 
1815, 2 octobre, 26 d\'ecembre.}

The manuscript of the prize was not published until 1827, when appeared the
first volume of \textit{M\'emoires des Savans \'etrangers} (Memoirs of
non-member scientists) printed since  the 1816  reorganization of the
Academy. It comprised a hefty  310 pages, including 189 pages with 20  
technical notes (the last 7 where added at various dates after 1815, but
all the material not contained in the 1815 manuscript is clearly
identified by the author).\footnote{The publication reference is found in our
reference list as Cauchy, 1815/1827, but what we discuss here is,
according to Cauchy, unchanged from the 1815 prized manuscript. On the circumstances of publication, see Belhoste, 1991; 
for the delay before publication, see also
Smithies, 1997:  \S\,2.1.}

It is interesting that Cauchy's prized manuscript  begins with a statement
of the problem to be solved:\\
{\footnotesize
\noindent
\textit{A massive fluid, initially at rest, and of an arbitrary depth, has been put in
motion by the effect of a given cause. One asks, after a determined time, the
shape of the external surface of the fluid and the speed of each molecule
located at the same surface.}\footnote{Original in French: Une masse fluide pesante, primitivement
  en repos, et d'une profondeur ind\'efinie, a \'et\'e mise en mouvement par
  le fait d'une cause donn\'ee. On demande, au bout d'un temps d\'etermin\'e,
  la forme de la surface ext\'erieure du fluide, et la vitesse de chacune
  des mol\'ecules situ\'ees \`a cette m\^eme surface.}  
}\\[-1.7ex] 

Although we have not
found this sentence in any of the minutes of the Academy or at its archives,
it is likely that it constitutes the Academy's detailed formulation of the
problem, which had been entrusted to Laplace.

Actually, Cauchy gave a more general treatment than had been  requested by the 
Academy, since he obtained results not only for the  \textit{surface} of the fluid, but 
also for  its \textit{ bulk}.  The memoir itself has three parts: the solution is in the third part, 
whereas the first two actually contain, in the intention of the Author,  
a sort of  preparatory 
material describing the initial state of the whole fluid and its later 
evolution.\footnote{For a short description of Cauchy's memoir,
see  Darrigol, 2005: 43--45. For a  detailed presentation of Cauchy's
paper in the context of 19th century work on wave propagation,
see Risser, 1925; Dalmedico, 1989.}

Here, it is the first section of the second part that interests us; it is  entitled \textit{On the equations that
subsist at any instant  of the motion for all the points within the mass
of the fluid.}\footnote{In modern scientific language:  The equations that
describe the fluid motion within its bulk at any time.} Our focus will
be entirely on this section and even more so on its very beginning, where
Cauchy obtains his  Lagrangian formulation of the 3D incompressible
Euler equations in terms of three invariants and uses it immediately to derive what is 
called Cauchy's vorticity formula relating the current and initial vorticity fields
through the Jacobian matrix of the Lagrangian map. In the next section
we turn to these matters.

\subsection{Cauchy's invariants equations and the Cauchy vorticity formula}
\label{ss:15-17}

\begin{figure}
\centerline{\includegraphics[height=2.0in]{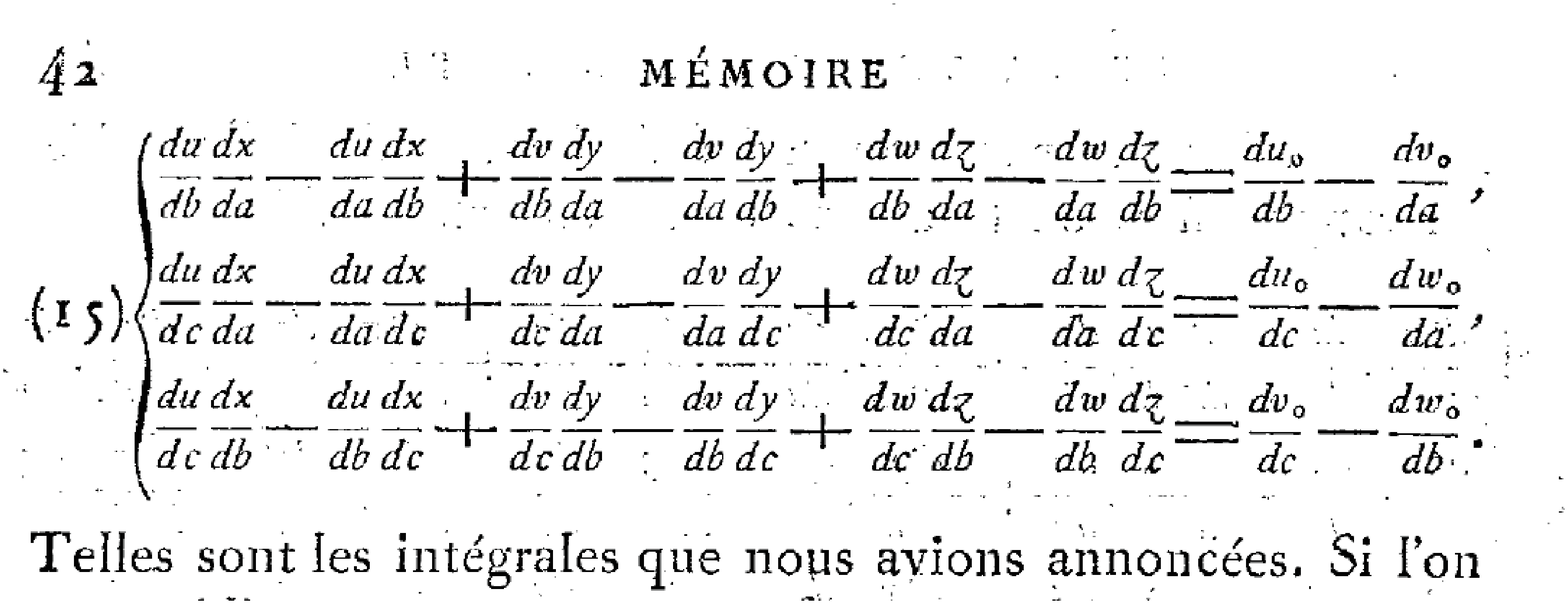}}
\caption{Cauchy's invariants equations (15) as they appeared in print in 1827.}
\label{f:cauchy15}\end{figure}

There is of course no real substitute for reading Section~I of the
second part of Cauchy's paper. It can be mostly understood without
prior reading of the first part. The notation is not too different from the 
modern one, except that, of course,  no vectors were used. For
illustration, Figure~\ref{f:cauchy15} gives Cauchy's key equation
(from the point of view of the present paper) as it was published in
1827. (Our attempts to retrieve the original hand-written manuscript of
1815 have failed.) In our description of the work we shall use modernized
notation.\footnote{Cauchy, 1815/1827: 35--49.}

Cauchy considers a 3D ideal incompressible fluid subject to an
external force. Analyzing the various forces acting on a ``molecule''
(i.e. a fluid particle), he derives his eq.~(4), which in our notation reads
\UB\partial_t\bv+(\bv\cdot\nabla)\bv=-\nabla p +\bF,\UE{EulerMom}
where, $\bv$ is the flow velocity, $\bF$ the external force, $p$ the
pressure (divided
by the constant fluid density, taken   here unity for convenience) and  $t$ the time. He then points out
that these equations coincide with those obtained by Lagrange by
another method. They coincide
also in their form and method of derivation with those obtained by
Euler and are presently called the Euler equations.\footnote{Euler, 1755; Lagrange, 1788: 453.}.

Cauchy, then, changes to Lagrangian variables, denoted by him
$(a,\,b,\,c,\,t)$ (here, $(\ba,\,t)$).  The Eulerian position $\bx$ becomes
then a function $\bx(\ba,t)$. Nowadays, the representation of the flow
in terms of the coordinates $(\bx,t)$ is called \textit{Eulerian} and,
when the coordinates $(\ba,\,t)$ are used, it is called
\textit{Lagrangian}; the (time-dependent) map $\ba \mapsto \bx$ is
called the \textit{Lagrangian map}.
The velocity and the acceleration of a
fluid particle are then $\dot \bx(\ba,\,t)$ and $\ddot\bx(\ba,\,t)$,
respectively, where the dot denotes the Lagrangian time derivative. The Euler equations
\eqref{EulerMom} state that the acceleration minus the external force
is balanced by minus the Eulerian gradient of the pressure. Making use
of the set of nine partial derivatives of the $\bx$s with respect to
the $\ba$s (now called the Jacobian matrix) Cauchy transforms the
Eulerian pressure gradient into a Lagrangian pressure gradient, here
denoted $\nL p$, and obtains his eq.~(7)
\UB\sum_{k=1}^3\left(\ddot{x}_k -F_k\right)\nL x_k=-\nL
p,\UE{diffWeber} where $x_k$ and $F_k$ denote the $k$th components of
$\bx$ and $\bF$, respectively. This is precisely Lagrange's
eq.~(D). Note that Lagrange first wrote the equations in Lagrangian
coordinates and then switched to Eulerian coordinates; Cauchy did it the
other way round.\footnote{Lagrange, 1788: 446.}

Then, Cauchy considers the condition of incompressibility, which he
first writes in Lagrangian coordinates. In modern terms, the Jacobian
of the Lagrangian map should be equal to unity for all $(\ba,t)$ (his
eq.~(9))
\UB\det\,(\nL\bx)=1.\UE{sdiff}
He also writes it in Eulerian coordinates (his eq.~(10))
\UB\nabla\cdot\bv=0,\UE{Eulerincomp}
an equation already found in Euler but which has been derived earlier
in the axisymmetric case by d'Alembert.\footnote{For a discussion
of Euler and d'Alembert's contribution to the incompressibility
condition,
see Darrigol \& Frisch, 2008: \S\S\,3, 4.}

Then, Cauchy observes that 
the two equations \eqref{diffWeber}-\eqref{sdiff} are not integrable,
but if one restricts oneself to external forces deriving from
a potential $\lambda$, then  \eqref{diffWeber} can be integrated once, as he will
show. Cauchy thus writes (his eq.~(11):
\UB\bF = \nabla \lambda.\UE{forcepotential}
He then rewrites \eqref{diffWeber} as (his eq.~(13))
\UB\sum_{k=1}^3\ddot{x}_k \nL x_k=-\nL
(p-\lambda).\UE{cauchy13}
Observe that the r.h.s is a Lagrangian gradient. Cauchy then applies
what we now call a (Lagrangian) curl to cancel out the r.h.s. He thus
obtains his eq.~(14)
\UB\nL \times\sum_{k=1}^3\ddot{x}_k \nL x_k=0.\UE{cauchy14}
He then notices that the three components of the l.h.s. of
\eqref{cauchy14}  are  exact time derivatives of three quantities
which thus must be time-independent. He easily identifies their
constant values to what we now call the initial vorticity $\bom_0 =
\nL\times \bv_0$. This way, Cauchy obtains his eq.~(15):
\UB\sum_{k=1}^3\nL\dot{x}_k\times\nL x_k=\sum_{k=1}^3\nL v_k\times\nL x_k=\bom_0,\UE{cauchy15}
and states: \textit{Telles sont les int\'egrales que nous avions
  annonc\'ees} (Such are the integrals that we had announced). Indeed, in
Lagrangian coordinates, the r.h.s. is time-independent. The constant
quantities in the l.h.s. are now usually called ``the Cauchy invariants,'' a
terminology we shall adopt. As we shall see in
Section~\ref{ss:hankel}, they are closely connected to the circulation
invariants of Helmholtz and Kelvin and are the three-dimensional generalization
of the two-dimensional vorticity invariant. As to the Cauchy equations \rf{sdiff}
  and \ref{cauchy15} for the Lagrangian map, which play a central role
  in the present paper, we shall refer to them as
``Cauchy's invariants equations'', to avoid any possible confusion with
``Cauchy's vorticity formula,'' discussed below. %Occasionally, we
  %shall also refer to the quantities that remain constant along fluid
  %particle trajectories as ``Cauchy invariants.''
 
Cauchy was obviously aware that he had succeeded in partially
integrating the equations of motion. However, modern concepts such as
invariants and their relation to symmetry/invariance groups would
emerge only about one century later (see
Section~\ref{ss:noether}). Nonetheless, for his invariants equations
Cauchy immediately found an application, which would become quite
famous (much more, so far, than Cauchy's invariants).

Starting from \rf{cauchy15}, written in terms of the velocity, Cauchy
reexpresses its Lagrangian space derivatives in terms of the Eulerian
ones
and the Jacobian matrix. He obtains for the l.h.s. of \rf{cauchy15} 
expressions which are linear in the components of the vorticity
$\bom$ (evaluated at the current time $t$) and quadratic in the
Jacobian matrix. He then solves these linear equations, using the fact
the Jacobian is unity. He thus obtains his eq.~(17):
\UB \bom = \bom_0\cdot \nL \bx, \UE{cauchy17}
or, with indices
\UB \bom_i = \sum_j\bom_{0j} \nL_j \bx_i. \UE{cauchy17i}
In modern terms this ``Cauchy vorticity formula'' states that the current
vorticity is obtained by multipying the initial vorticity by the
Jacobian matrix. Cauchy gives a rather low-key application of his
formula, which is consistent with the context of the prize: in the
first part of his memoir, he had envisaged a mechanism of setting the
fluid
in motion impulsively that would produce a flow initially potential
and thus with no vorticity. His formula then implies that the flow
would have no vorticity at any instant of time. In the language
of the time this was expressed by stating that $\bv\cdot {\rm d} \bx$
is a ``complete differential.''

Nonetheless, Cauchy was certainly aware of Lagrange's theorem, which
states  that an ideal flow initially potential, stays potential at
later times. Lagrange's  proof used Eulerian coordinates and assumed that the
velocity
could be Taylor-expanded in time to arbitrary orders. Lagrange then
showed that if the vorticity vanishes initially, so will its
time-Taylor coefficients of arbitrary orders.\footnote{Lagrange, 1788: 458--463.}

Cauchy's proof requires only a limited smoothness of the flow (he does
not state how much) and it must have appeared to the readers at the
time that, as long as the Jacobian matrix exists, then the persistence
of potential flow will hold.  Stokes observed that the vanishing of
all the Taylor coefficients does not imply the vanishing of a function
(giving well-known examples such as $\ue ^{-1/t^2}$ near $t=0$); he
thus considered Cauchy's proof more general than that of
Lagrange. Today, we know that a flow with an initial velocity field
that is ``moderately smooth in space'' (just a little more regular
than once differentiable in space), will stay so for at least a finite
time, during which its temporal smoothness in time \textit{in Eulerian
  coordinates} is not better than its spatial smoothness, rendering
Lagrange's argument inapplicable, whereas Cauchy's proof only requires spatial
differentiability of the Lagrangian map, which is now known to hold
with a moderately smooth initial velocity field.\footnote{Stokes, 1845:
  305. On spatial differentiability see, e.g. Zheligovsky \&
  Frisch, 2014 and references therein.}

One reason why Cauchy's vorticity formula \rf{cauchy17} is very well known today
is that it applies not only to the vorticity in an ideal fluid,
but also to a magnetic field in ideal conducting fluid flow governed
by the magnetohydrodynamic (MHD) equations. In modern mathematical
language, both the vorticity and the magnetic field are transported
2-forms. In an ordinary fluid, 
one cannot prescribe the velocity and the
vorticity independently, but in a conducting fluid one can prescribe
the velocity and the initial magnetic field independently in the limit
of weak magnetic fields, when studying the kinematic dynamo problem.\footnote{On MHD and dynamo theory, see Moffatt, 1978.}

All this explains that there has been a strong interest, particulary
in recent years, in Cauchy's vorticity formula \rf{cauchy17}. Our main
focus in this paper are the Cauchy invariants equations \rf{cauchy15}.
One cannot describe the history of the Cauchy's vorticity formula
without mentioning his invariants equations, and thus we cannot
completely disentangle the histories of their citations. However, 
nowadays, most derivations of Cauchy's vorticity formula use a
much shorter route, based on the Eulerian vorticity equation of Helmholtz
\UB\partial_t \bom +\bv\cdot\nabla \bom =\bom\cdot \nabla\bv,\UE{helm} and thus bypass
Cauchy's invariants equations.\footnote{Helmholtz, 1858.}

The rest of Cauchy's Section I of Part II (pp.~44-49)  is
devoted to the  case of potential flow and does not concern us here.

\section{19th century}
\label{s:19th}

\subsection{Hermann Hankel}
\label{ss:hankel}

\begin{figure}
\centerline{\includegraphics[height=5in]{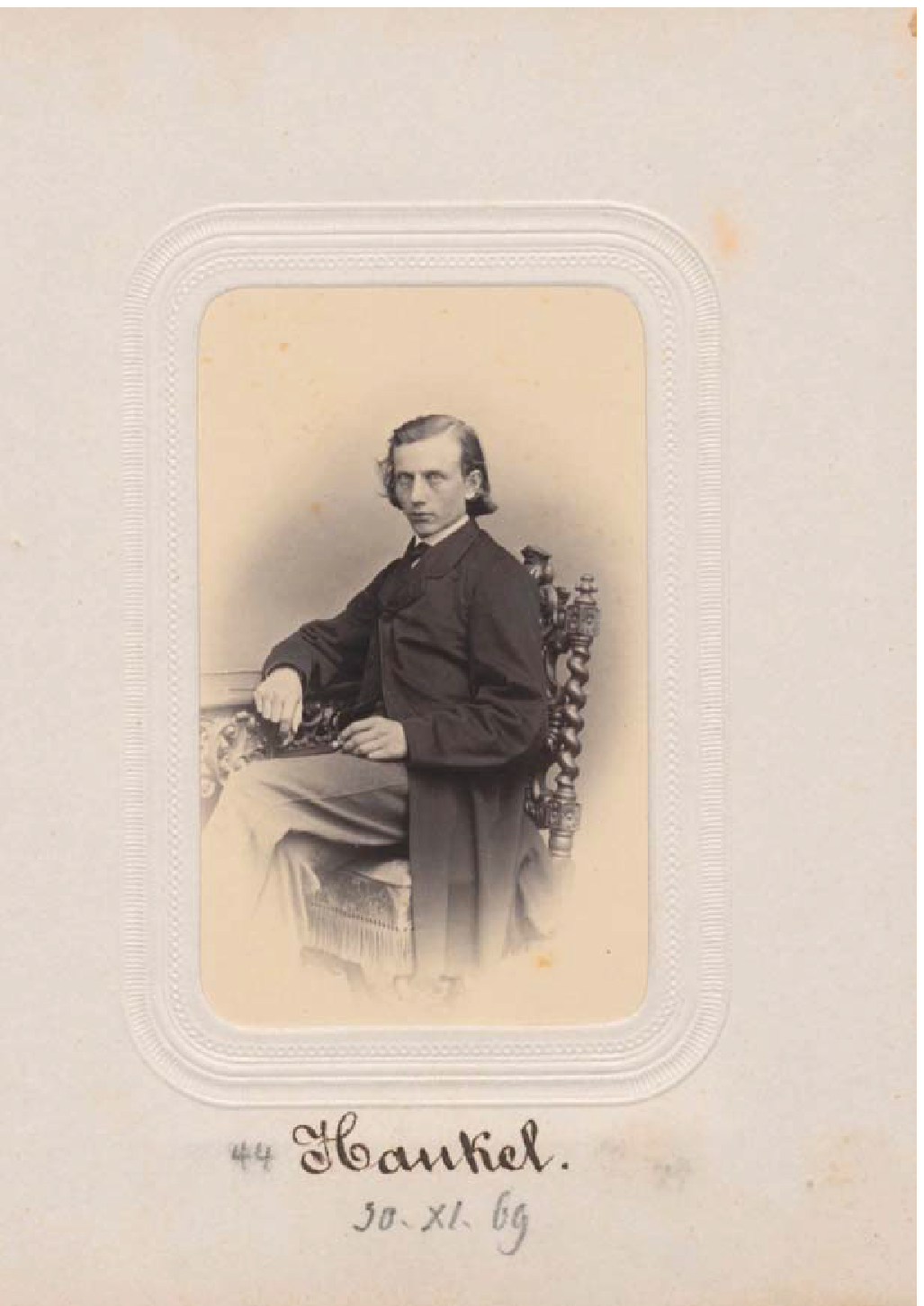}}
\caption{Hermann Hankel, Tobias-Bild Universit\"atsbibliothek
  T\"ubingen. By permission.}
\label{f:hankel}\end{figure}

Hermann Hankel (1839--1873) was a German mathematician who studied with
August Ferdinand Moebius, Bernhard Riemann, Karl Weierstrass and Leopold Kronecker. Our sources on his life
are the obituary by Wilhelm von Zahn, a 19th century biography by
Moritz Cantor, a 20th century short biography  by  Michael Crowe and an assessment of his mathematical contribution
from a modern perspective by Antonie Frans Monna.\footnote{von Zahn,
  1874. Cantor, 1879. 
 Crowe, 2008. Monna, 1973.}

Hankel is quite well-known for work on Hankel matrices, on
Hankel transforms and on Hankel functions. He was in addition also involved very seriously in the
history of mathematics and has left a book \textit{Zur Geschichte der
  Mathematik in Alterthum und Mittelalter} (On the History of
Mathematics in the Antiquity and the Middle Ages) which includes,
among other things,  one of the
first studies bringing out the major contributions of Indian
mathematics.\footnote{Hankel, 1874.} What interests us here is Hankel's work on fluid
dynamics contained in a  manuscript prized by G\"ottingen University \textit{Zur allgemeinen
  Theorie der Bewegung der Fl\"ussigkeiten} (On the general theory of
motion of fluids).

To understand the context of this prize, let us recall that in 1858 Hermann
Helmholtz (1821--1894) wrote a major paper about vortex motion in three-dimensional
incompressible ideal flow that generated considerable
interest.\footnote{Helmholtz, 1858; for his work on vortex dynamics,
  see also the biography by K\"onigsberger, 1902--1903, 
  English version 1906: 167--171 and Darrigol, 2005:~Chap.~4.
For an extensive bibliography of vortex flow since Helmholtz, see
Meleshko \& Aref, 2007.}

A key result of his work, stated in modern language, is that the flux
of the vorticity through an infinitesimal piece of surface is a
Lagrangian invariant. This result,  known as Helmholtz's second
theorem, 
immediately implies, by adding up infinitesimal
surface elements, that the same holds for a finite surface;
moreover, by using Stokes's theorem one obtains Kelvin's
circulation theorem. Helmholtz's derivation of his result is to a
large extent resting on an Eulerian approach and begins with the
establishment of the aforementioned Eulerian vorticity
equation. Furthermore, Helmholtz's derivation, mostly written for
physicists, was a bit heuristic. %It is thus not very surprising

On June 4, 1860 G\"ottingen University  (Philosophische Facult\"at
der Georgia Augusta) set up a prize, intended to stimulate interest
in Lagrangian approaches and in particular to give such a derivation of 
Helmholtz's invariants :\\
{\footnotesize
\noindent
\textit{The general equations for determining fluids motions may be  given  
in two ways, one of  which  is Eulerian, the other one is Lagrangian. The illustrious  Dirichlet pointed out in the posthumous unpublished  paper
``On a problem of hydrodynamics'' the until now almost completely overlooked advantages  of the
Lagrangian way, but   he was  prevented from unfolding this way further  by a
deadly illness. So, this institution asks for a theory of fluid motion  based on
the equations of Lagrange, yielding, at least, the laws of vortex motion 
already derived in another way by the illustrious Helmholtz.}\footnote{Aequationes generales 
motui fluidorum determinando inservientes
duobus modis exhiberi possunt, quorum alter Eulero, alter Lagrangio
debetur. Lagrangiani modi utilitates adhuc fere penitus neglecti 
clarissimus Dirichlet indicavit in commentatione postuma ``de problemate
quodam hydrodynamico'' inscripta; sed ab explicatione earum uberiore
morbo supremo impeditus esse videtur. Itaque postulat ordo theoriam
motus fluidorum aequationibus Lagrangianis superstructam eamque eo
saltem perductam, ut leges motus rotatorii a c1arissimo Helmholtz
alio modo erutae inde redundent.}  
}\\[-1.7ex] 

The prominent reference to Dirichlet can be understood as follows.
Johann Peter Gustav Lejeune Dirichlet (1805--1859) was an important German 
mathematician who, from 1855 to his
death,  succeeded Carl Friedrich Gauss in G\"ottingen. He had also a strong interest in hydrodynamics. In 1856--1857 he wrote
an unfinished paper ``Untersuchungen \"uber ein Problem der Hydrodynamik''
(Investigation of a problem in hydrodynamics).  Dirichlet asked
the German mathematician
Richard Dedekind (1831--1916), at that time professor in G\"ottingen,  
 to help him with the work, but was not able to 
finish completely before he died. Major
pieces of the work were found by Dedekind who published them in 1859
with the above title, followed by ``From his legacy, edited by
R. Dedekind.'' In the Introduction  Dirichlet pointed out 
that Lagrange himself was not too keen to advocate  the use of what was
later called Lagrangian coordinates, found by Lagrange to be a bit complicated.
Eulerian coordinates quickly grew in favour. Dirichlet, however
observed that Eulerian coordinates have their own drawbacks, particularly
when the volume occupied by the fluid changes in
time.\footnote{Dirichlet, 1859. On the role of Euler in introducing
  both the ``Eulerian'' and the ``Lagrangian'' coordinates, see Hankel,
  1861: 3 and Truesdell, 1954a: 30 (footnote 2).}

We now turn to Hankel's prize-winning manuscript. It carried the epigraph \textit{Tanto utiliores sunt notae, 
quanto magis exprimunt rerum relationes} (The more signs express relations
among things, the more useful they are).
The manuscript was written in Latin, but Hankel got the permission to print a slightly
edited German translation, which was also published as a book. In 1863 a four-page review was published 
in \textit{Fortschritte der Physik} (Progress in Physics); it was signed,
anonymously as ``Hl.'' and will be cited
below as \textit{Fortschritte}.\footnote{Hankel, 1861. Anonymous, 1863.}

This work is of
particular
interest, not only because it is the first time that it was
shown that the Helmholtz invariants are directly connected
with the Cauchy invariants, but also because it gives the first
derivation
of what is generally called the Kelvin circulation theorem. Hankel discusses both compressible and
incompressible fluids, but here it suffices to consider the latter.

To derive the Helmholtz invariants, Hankel first establishes 
Cauchy's invariants equations \eqref{cauchy15}, for which he refers to Cauchy's 1815/1827
prized paper. Hankel then rewrites this, in modernised notation,
\UB\nL \times \sum_{k=1}^3\dot{x}_k\nL x_k=\bom_0,\UE{hankelst}
where $v_k = \dot{x}_k$ are the components of the velocity $\bv$. This
is eq.~(3) of Hankel's \S\,6. The equation is not to be found in Cauchy's 1815/1827
paper but, given that Cauchy obtained his invariants equations (here, \rf{cauchy15}) by taking the
Lagrangian curl of \rf{diffWeber} and then integrating over the time
variable, it is not really surprising that the left-hand-side of
\rf{cauchy15} is a  Lagrangian curl.

Hankel's next step is to consider in Lagrangian space a connected
(\textit{zusammenh\"angende}) surface, here denoted by $S_0$, whose
boundary is a curve, here denoted by $C_0$, and to apply
to \rf{hankelst} what was later to be called the Stokes theorem,
relating the flux of the curl of a vector field across a surface
to the circulation of the vector field along the boundary of this surface. Hankel
could not easily be aware of what Thomson and Stokes had done before
on the subject and of the much earlier work of Ostrogradski and thus he devotes his
\S\,7 to proving the Stokes theorem.\footnote{Concerning the history of the
  Stokes theorem, see Katz, 1979.}

In \S\,8, Hankel then infers that the flux through $S_0$ of the 
r.h.s. of \rf{hankelst}, namely the initial vorticity,
is given by the circulation along $C_0$ of the initial velocity $\bv_0$:
\UB \int_{C_0}\bv_0\cdot d\ba=\int_{S_0}\bom_0
\cdot \bn_0 \,d\sigma_0,\UE{hankelhtu}
where $\bn_0$ denotes the local unit  normal to $S_0$ and $d\sigma_0$ the surface
element. This is the unnumbered equation near the top of his p.~38.
 Then, in \S\,9, he similarly handles the l.h.s. of \rf{hankelst}
and first notices that 
\UB\sum_{k=1}^3 v_k\nL x_k \cdot d\ba= \bv\cdot d\bx.\UE{hankelndt}
This is the third equation before eq.~(2) of his \S\,9. He thus
obtains the Eulerian circulation, an integral over the curve $C$ where are
located at the present time the fluid particles initially on $C_0$:
\UB \int_C \bv\cdot d\bx = \int_{S_0}\bom_0
\cdot \bn_0 \,d\sigma_0,\UE{hankelndu}
which is the unnumbered equation just before eq.~(2) of his \S\,9.
Eq.~\rf{hankelndu}, together with \rf{hankelhtu} is clearly the
standard circulation theorem, generally associated to the name of Kelvin.

Hankel does not seem to wish highlighting this result, hence the
unnumbered equations. Also, the result is not mentionned in the
\textit{Fortschritte} review.  Nevertheless, the fact that Hankel proved the
circulation theorem eight years before Kelvin did not escape the attention
of Truesdell, who even proposed calling it the ``Hankel--Kelvin
circulation theorem.'' However, Truesdell did not explain how Hankel
proceeded and furthermore never cited Cauchy's invariants equations,
but just the Cauchy's vorticity formula. This could be the reason why Truesdell's
rather justified suggestion did not seem to have many followers, one exception
being a book on ship propellers by Breslin and
Andersen who were aware of Truesdell's suggestion.  
One further
application by Hankel of the Stokes theorem, gives him the constancy
in time of the flux of the vorticity through any finite surface moving
with the fluid. Finally, he lets this surface shrink to an
infinitesimal element and obtains Helmholtz's theorem.\footnote{On
  Thomson's (Lord Kelvin) own derivation of the circulation theorem,
  see Thomson (Lord Kelvin), 1869. On the naming as ``Hankel--Kelvin
  theorem'', see Truesdell, 1954a:
  \S\,49. Breslin \& Andersen, 1996: 497.}

To conclude this section on Hankel, we ask: how much was his work on 
hydrodynamics remembered?  An interesting case is that of Heinrich Martin Weber
(1842--1913), who was quite close to  Riemann. In 1868, Weber wrote
a paper titled \textit{\"Uber eine Transformation der hydrodynamischen Gleichungen}
(On  a transformation of the equations of hydrodynamics) which, from the point
of view of its scientific content, is very closely related to Cauchy's 1815
invariants equations and even more so to Hankel's 1861 reformulation 
\rf{hankelst}.  Specifically, by ``decurling'' 
\rf{hankelst} in Lagrangian coordinates, one obtains
\UB\sum_{k=1}^3\dot{x}_k\nL x_k=\bv_0 -\nL W,\UE{weber1}
where $\bv_0$ is the initial velocity and $W$ a scalar function, here called ``the Weber function''. Actually,
Weber showed that $W$ is the time integral from $0$ to $t$, in Lagrangian 
coordinates, of $p-(1/2)|\bv|^2$, where $p$ is the pressure and $\bv$ the
velocity. Weber derived his equation by a clever transformation
of Lagrange's equation \rf{cauchy14}, now called the  ``Weber transform.''
Weber did cite Hankel but without an actual reference and just as a person who
had pointed out that the so-called Eulerian and Lagrangian coordinates
both were first introduced by Euler (a statement made by Hankel but attributed
by him to his advisor, Riemann).\footnote{Weber, 1868.}

\begin{figure}
\centerline{\includegraphics[height=4in]{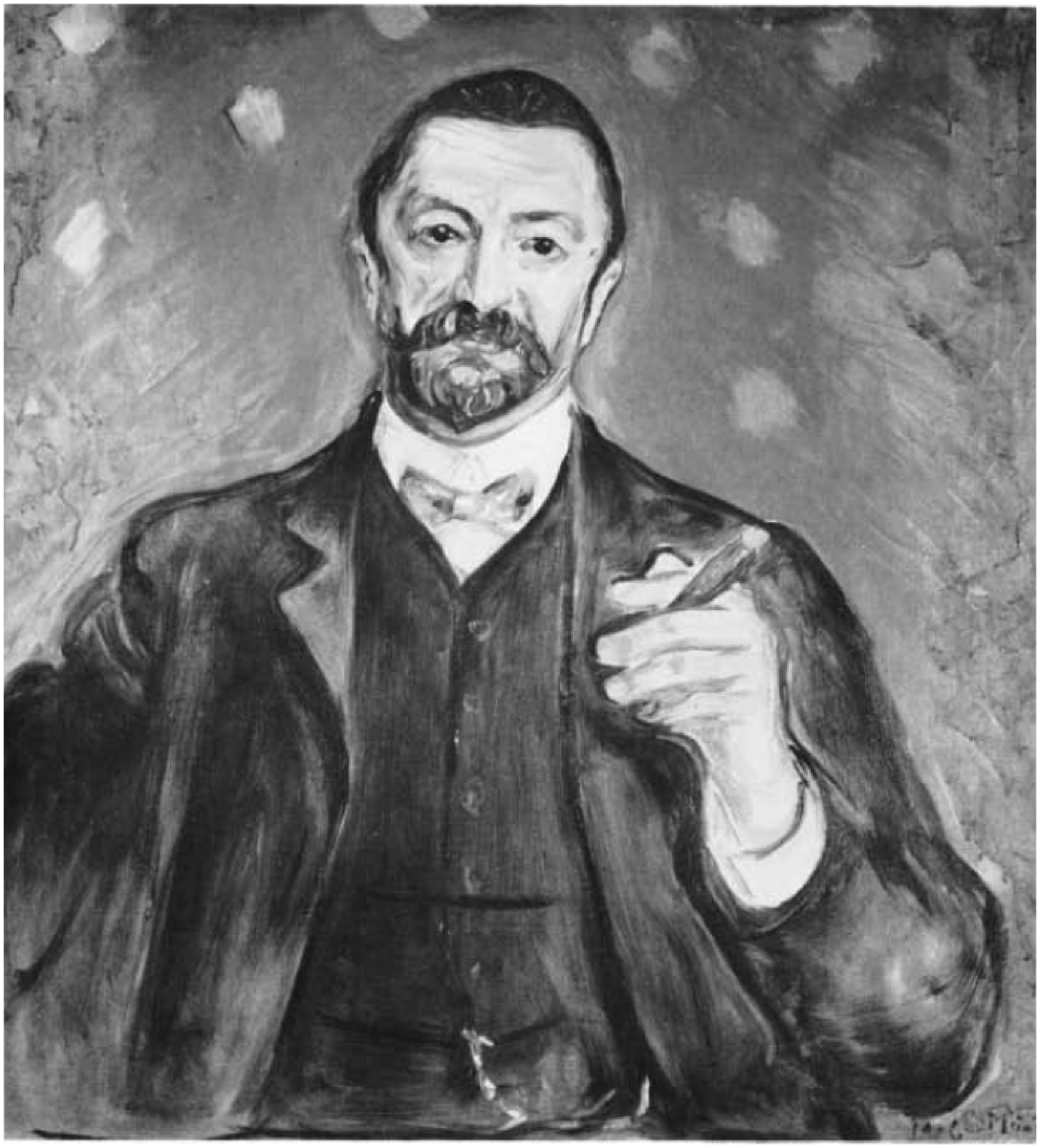}}
\caption{Portrait of Felix Auerbach by Edvard Munch (1906). Private collection. Reproduced from the Burlington Magazine (1964) by permission.}
\label{fig3}\end{figure}
Felix Auerbach (1856--1933),\footnote{He and his wife took their lives
on February 26, 1933.} a German scientist with wide-ranging interests in all areas of physics,
in hydrodynamics, in architecture and painting,  wrote in his early days
\textit{Die theoretische Hydrodynamik
nach dem Gange ihrer Entwickelung in der neuesten Zeit, in K\"urze
dargestellt} (A brief presentation of theoretical hydrodynamics, following its evolution
in the most recent times). This manuscript won the Querini Stampalia Foundation
prize of the Royal Venetian Institute of Sciences, Letters and Arts
(\textit{Atti del Reale Istituto Veneto di Scienze, Lettere ed Arti}) on the assigned theme
of ``Essential progress of theoretical hydrodynamics.''  Several pages are devoted to Hankel's
hydrodynamics work and the manuscript  contains also a brief reference
to Cauchy on p.~34.\footnote{Auerbach, 1881.}

\section{Stokes, L\'evy, Lamb, Andrade and Appell}
\label{ss:stokesetal}

\begin{figure}
\centerline{\includegraphics[height=4in]{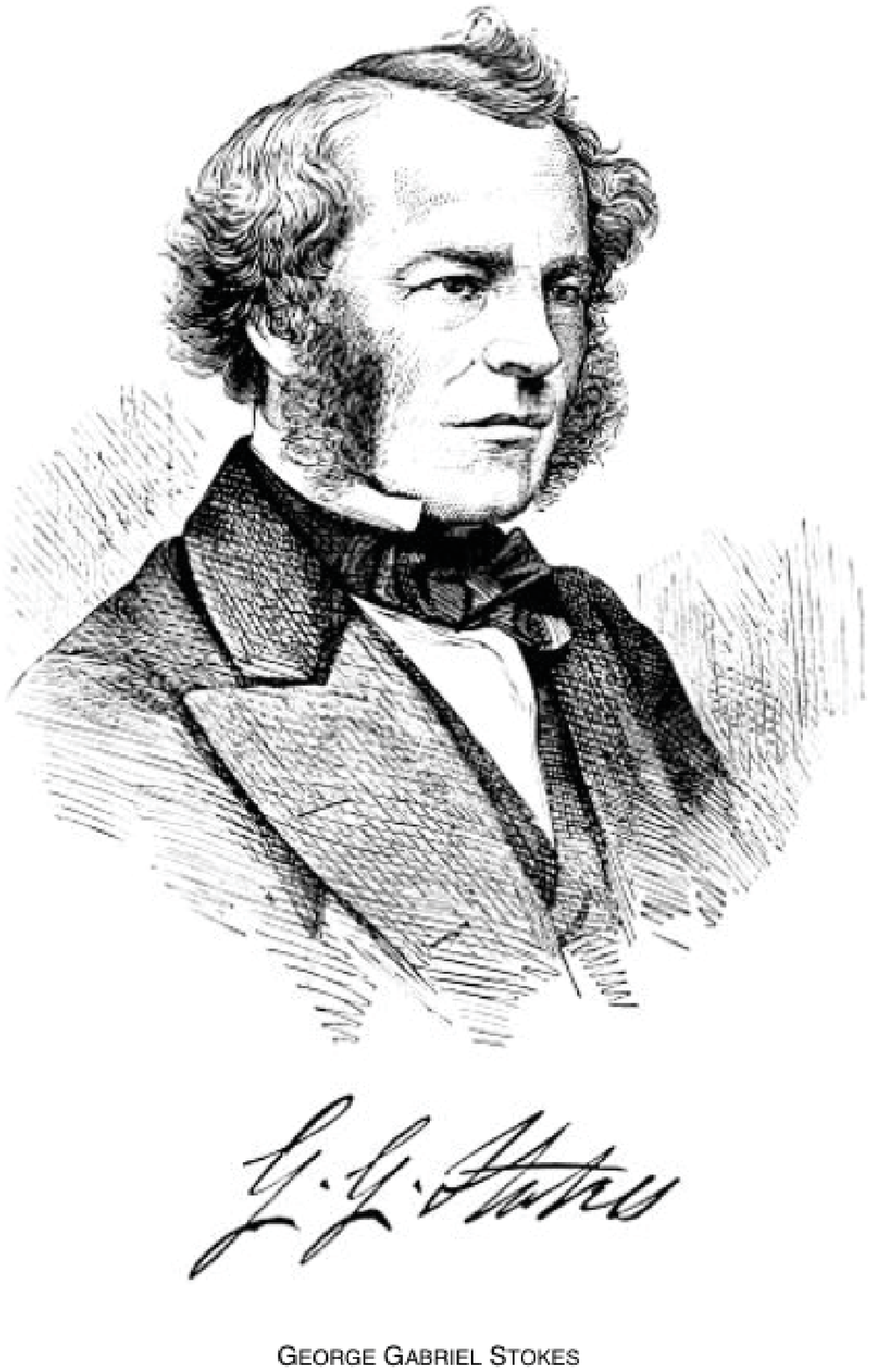}}
\caption{From Wikimedia Commons (Public domain).}
\label{f:stokes}\end{figure}
An Irish 
physicist and mathematician, George Gabriel Stokes (1819--1903) 
spent all his career at the University of Cambridge in England and was
considered a leading British scientist, particularly so for many
contributions to the dynamics of both ideal and viscous
fluids.\footnote{For the biography of Stokes, see Parkinson, 2008,
  Wood, 2003,
Wilson's 2011 two volumes  on the correspondence between Stokes and Kelvin with
many biographical elements (Stokes \& Thomson, 1846--1869,
1870--1903), and the obituary by
Rayleigh (Strutt (Lord Rayleigh), 1904).}

Stokes followed rather closely the work of French mathematicians and
physicists and made genuine efforts to cite other scientists's work. To the best of our knowledge, he was the first to realize
the importance of the discovery of the Cauchy invariants 
(called by Stokes ``integrals'') and of the ensuing Cauchy formula
\rf{cauchy17} for the vorticity. In three papers in the late 1840s,
Stokes discussed various proofs of  Lagrange's theorem on the
persistence in time of potentiality for 3D incompressible
flow.\footnote{Stokes, 1845, 1846, 1848.}

In particular, in the 1848 paper ``Notes on Hydrodynamics IV'' Stokes
described in detail Cauchy's proof of Lagrange's theorem and also gave an
alternative proof of his own. Concerning Cauchy's proof, Stokes
wrote:\\ {\footnotesize
\noindent
\textit{The theorem considered follows as a particular consequence from M. Cauchy's integrals. As however the equations
employed in obtaining these integrals are rather long, and the
integrals themselves do not seem to lead to any result of much
interest except the theorem enunciated at the beginning of this
article.}
}\\%[-1.7ex]
Stokes also gave an alternative proof of his own, not using the Cauchy 
integrals.

However, as observed by Meleshko and Aref, in 1883 when Stokes edited his ``Mathematical and Physical Papers'',
he introduced a footnote, refering to an added note at the end of the paper.\footnote{Stokes, 1883. Meleshko \& Aref, 2007: \S\,2.}
This note begins as follows:\\
{\footnotesize
\noindent
\textit{It may be noticed that two of Helmholtz's fundamental propositions 
respecting vortex motion follow immediately from
Cauchy's integrals; or rather, two propositions the same as those
of Helmholtz merely generalized so as to include elastic fluids
follow from Cauchy's equations similarly generalized.}
}\\[-1.7ex]

The two propositions are (i) that ``the same
loci of particles which at one moment are vortex lines remain
vortex lines throughout the motion'' and (ii) in modernised language, that the product of 
the modulus of the vorticity and of the area of a
perpendicular section of an infinitesimal vortex tube does not change
in time while following the Lagrangian motion.
Actually Stokes's statement should not be misread: he mentions
``Cauchy's integrals'' but, in 1883, Stokes understands by this only the Cauchy vorticity 
formula  \rf{cauchy17}, which of course was derived from Cauchy's
invariants
equations \rf{cauchy15}.

Maurice L\'evy (1838--1910) was a French engineer and specialist of
continuum mechanics. In 1890 he gave a 
lecture on ``Modern hydrodynamics and the hypothesis of action at a
distance''
at the Coll\`ege de France where he was a professor.\footnote{L\'evy, 1890.} On the first page
of the published version, L\'evy writes:\\
{\footnotesize
\textit{\noindent The admirable properties of vortices were discovered only in
1858 by Helmholtz, although they merely express the intermediate
integrals of Lagrange's hydrodynamical equations, discovered by Cauchy...}
}\footnote{Les admirables propri\'et\'es des tourbillons
n'ont \'et\'e d\'ecouvertes qu'en 1858 par Helmholtz, bien
qu'elles n'expriment pas autre chose que les int\'egrales  interm\'ediaires
des \'equations de l'Hydrodynamique de Lagrange,  d\'ecouvertes par
Cauchy...}\\[-2.5ex]

L\'evy observed that Cauchy wrote \textit{three} 
(scalar)
conservation laws; together with the condition of incompressibility
this makes \textit{four} equations for the three components of the
Lagrangian map. L\'evy stated that the equations are actually
compatible (this follows from the fact that the Lagrangian
divergence of Cauchy's three invariants vanishes). It is of particular
interest that L\'evy very much highlighted what we would today call
the \textit{nonlocal} character of the equations of incompressible
fluid dynamics. This did not seem to him in violation of any
known mechanical principle. Of course, such observations were made fifteen
years before the birth of relativity theory.  Today, we know that 
the nonlocal character stems from taking the limit of vanishing
Mach number for a slightly compressible fluid, a limit that amounts
to letting the speed of sound tend to infinity.

Horace Lamb (1849--1934), a British applied mathematician, wrote one
of the most authoritative treatise on hydrodynamics, with editions
ranging from 1879 to 1932 (the title ``Hydrodynamics'' was used only
from 1895).   From 1895, Lamb had Cauchy's invariants equations 
\rf{cauchy15}  but
only as an intermediate step to obtain Cauchy's vorticity formula \rf{cauchy17} and Cauchy's
derivation of Lagrange's  theorem. There is no
mention of ``integrals'' or ``invariants'', although mere inspection
of the equations makes it clear that we are here dealing with integrals
of motion, as Cauchy himself pointed out in 1815.\footnote{Lamb, 1932:
  \S\,146, p.~204 first unnumbered equation; also \S\,143, p.~225 first unnumbered equation of the 1895
  edition.}

Paul \'Emile Appell (1855--1930), a French mathematician with a talent
for simple and illuminating writing, published in 1897 a paper where
he gave an elementary and immediate interpretation of Cauchy's
equations,
leading to the fundamental theorems of the theory of vortices. He
first
rederived Cauchy's invariants equations \rf{cauchy15}. He then considered the
following first-order differential form
\UB\bv\cdot d\bx - \bv_0\cdot d\ba,\UE{appell-poincare}
where $\bv_0$ is the initial velocity, and showed that as a consequence
of \rf{cauchy15} it is an exact differential $dV$ of some function
$V$. (Actually $V=-W$ where $W$ is the Weber function defined near the end of 
Section~\ref{ss:hankel}.) Since the integral of such an exact form on a closed contour
vanishes
with suitable connectedness and regularity assumptions, Appell
immediately
obtained the Hankel--Kelvin circulation theorem. Appell pointed out
that here he was  just following Poincar\'e's \textit{Th\'eorie des tourbillons} (Lectures on
vortices). Actually, Poincar\'e's derivation is a bit more
mathematical 
and quite close to Hankel's derivation of the circulation theorem
(see Section~\ref{ss:hankel}). 
The derivation is again based on the Cauchy invariants equations
\rf{cauchy15} for which Poincar\'e cites Kirchhoff's ``Lectures on mathematical
physics (mechanics)''. The latter writes indeed  Cauchy's invariants equations
\rf{cauchy15} [Lecture 15, \S\.3,  eq.~(14) on p.~165]  but does not
give any reference.\footnote{Appell, 1897. Poincar\'e, 1893. Kirchhoff, 1876.}

Jules Andrade (1857--1933), a French specialist of mechanics and
chronometry, published in 1898 ``\textit{Le\c{c}ons de m\'ecanique physique}''
(Lectures on physical mechanics). Its Chapter VI was devoted to fluid
dynamics. On p.~242, Andrade derived the Cauchy invariants equations
\rf{cauchy15}, which he called  ``Cauchy's intermediate integrals''.
Andrade  also derived Cauchy's vorticity formula \rf{cauchy17} and Lagrange's theorem, closely 
following Cauchy. Andrade then stated \textit{Les th\'eor\`emes d'Helmholtz sont aussi
referm\'es dans ces \'equations, mais Cauchy ne les a pas
aper\c{c}us.} (Helmholtz's theorems are also contained in these
equations but Cauchy did not perceive them). He then showed how to derive Helmholtz's
result along more or less the lines used by Stokes in his 1883 added note
(see above).\footnote{Andrade, 1898.}

\section{20th century}
\label{s:20th}

The material is here separated into three subsections: Section~\ref{ss:noncauchy}
has not only Cauchy's work on the invariants equations forgotten,
but the invariants themselves never mentioned. In
Section~\ref{ss:noether}, we find the independent rediscovery of the
Cauchy invariants by application of Noether's theorem. Eventually, in 
Section~\ref{ss:rebirth} everything
will reconnect in the late 20th century.

\subsection{The (Cauchy) invariants fallen into oblivion}
\label{ss:noncauchy}

 The 20th century was to see a tremendous rise of research in fluid dynamics,
 driven to a significant part by the needs of the blossoming aeronautical
 industry. For this, the study of flow constrained by external or internal
 boundaries -- with viscous boundary layers of the kind introduced by 
Prandtl in 1905 -- was essential.\footnote{Prandtl, 1905. For a review of German fluid mechanics at the dawn of the 20th century,
see Eckert, 2006.}
%Prandtl, L. 1905. Uber Fl\"ussigkeitsbewegung bei sehr kleiner Reibung. 
%In Ver-handlungen des III. Internationalen Mathematiker-Kongresses, 
%Heidelberg, 484--491. Also in Selected Works, Walter Tollmien, 
%Hermann Schlichting and Henry G\"ortler, eds. vol. 2,  575--584.
%Eckert, M. 2006 The Dawn of Fuid Dynamics. A Discipline between Science
%and Technology. Wiley-VCH Verlag (Weinheim).
 
Inclusion of viscous effects requires the use of the
 Navier--Stokes equations, which are somewhat easier to study in Eulerian
 coordinates.  Mathematical issues, relating to ideal and viscous fluid flow,
 such as the well-posedness of the fluid dynamical equations, started being
 addressed with the new tools of functional analysis.  For example,
 Lichtenstein gave the first proof of the well-posedeness for at least
 a finite time of the three-dimensional incompressible Euler equations with
 sufficiently smooth initial data. Then, H\"older and Wolibner 
independently showed
 that, under suitable conditions, the two-dimensional incompressible Euler
 equations constitute a well-posed problem for all times. Leray obtained
 similar results in the viscous case and introduced the important concept of
 ``weak solutions'' which need not be differentiable.\footnote{Lichtenstein, 
1927. H\"older, 1933. Wolibner, 1933. Leray, 1934.}

We do not give here other details.\footnote{See, e.g., Rose \& Sulem, 1978; Majda \& Bertozzi, 
2002.} We stress that these results were generally obtained using
 Eulerian coordinates, with an occasional excursion into Lagrangian
coordinates by Lichtenstein. It seems that during the 20th
 century Cauchy's invariants equations \rf{cauchy15} were hardly
 used for mathematical studies. On the one hand, this could be because of the general
belief at that time, going back to Lagrange, that (what we now call) Lagrangian
coordinates are unnecessarily complicated; as we already stated, 
the questionable character of this belief was underlined by 
Dirichlet.\footnote{Dirichlet, 1859.} On the other hand, it may be that
Cauchy's Lagrangian formulation through \rf{sdiff}-\rf{cauchy15} just drifted 
into oblivion.

In order to understand better what had happened, we examined, in addition
to the mathematical papers already cited,  a considerable number of 
major fluid
mechanics  textbooks published in the 20th century and looked for citations
of Cauchy's 1815/1827 work. A fully relevant citation would not only
have Cauchy's invariants equations \rf{cauchy15}, but also stress, as 
Cauchy did, that they define Lagrangian invariants. Here a word of caution
is required: since Stokes in 1883, several authors have referred
to the Cauchy vorticity formula \rf{cauchy17} as ``Cauchy's integrals.'' Cauchy
used the word ``integrals'' in connection with \rf{cauchy15}
and not \rf{cauchy17}.   We failed to find
any truly relevant citations before the very end of the 20th century,
although Cauchy's invariants equations (or a 2D instance) were
rediscovered independently. Hereafter,
we indicate some of the partially relevant findings.

 Lamb's treatment of Cauchy remained 
exactly what it was in 1895 (see Section~\ref{ss:stokesetal}) with little
emphasis on \rf{cauchy15}.

Lichtenstein, in addition to his pioneering papers on
the mathematical theory of ideal flow, published several books. In
1929 he produced volume XXX \textit{Grundlagen der Hydrodynamik}
(Foundations of hydrodynamics) of an encyclopedia of mathematics with
emphasis on applications, edited by Richard Courant. Here, in Chapter
10, Cauchy's invariants equations appear briefly [his eq.~(54) for the
  case $\rho=1$] but are not directly attributed to Cauchy and no use
is made of the invariance other than deriving the Cauchy vorticity
formula.\footnote{Lichtenstein, 1929: 397.}

 Sommerfeld's 1945  ``Mechanics of deformable bodies'' and
Landau and Lifshitz's 1944 first edition of ``Fluid Mechanics'' seem to contain neither
the Cauchy invariants nor the Cauchy vorticity formula. In his Encyclopedia of
Physics article on the foundations of fluid dynamics, 
Oswatitsch  followed
rather closely Cauchy's original derivation of the vorticity formula
but  skipped the
Cauchy invariants (a similar treatment with some allowance for turbulent
fluctuations is made by Goldstein). We have already mentioned in Section~\ref{ss:hankel} that Truesdell cited Hankel quite extensively, but cited Cauchy only for
his vorticity formula. Ramsey has \rf{cauchy15} but, again, only
as an intermediate step in proving Cauchy's vorticity
formula. Batchelor in his ``Introduction to Fluid Dynamics'' derived the vorticity formula and attributed it
to Cauchy. Finally, Stuart and Tabor in their introductory paper to the Theme Issue of 
\textit{Philosophical Transactions}, devoted to the Lagrangian description of fluid motions, have Cauchy's invariants equations: these are their
equations (2.13)--(2.15), which were here derived from the Cauchy
vorticity formula; this amounts to retracing Cauchy's steps in reverse.
It is not mentioned that the resulting equations were already in 
Cauchy 1815/1827.\footnote{Sommerfeld, 1950. Landau \& Lifshitz, 1959. 
Oswatitsch, 1959. Goldstein, 1938: Chap.~5, \S\,84. Truesdell, 1954a. Ramsey,
1913: Chap.~2, 22. Batchelor, 1967: 276. Stuart \& Tabor, 1990.}

\subsection{The (Cauchy) invariants rediscovered and Noether's theorem}
\label{ss:noether}

As we shall see, the rebirth of Cauchy's invariants at the very end
of the 20th century and the beginning of the 21st
(Section~\ref{ss:rebirth}) was preceded by rediscoveries (without
Cauchy being named). The most widely known  of these rediscoveries used a novel tool developed in the early 
20th century. Emmy Noether (1882--1935) is recognized  as one of
the most important mathematicians of all times for her work in
algebra.\footnote{See, e.g., Dick, 1970 and the Wikipedia article 
at \url{http://en.wikipedia.org/wiki/Emmy_noether}.}

\begin{figure}
\centerline{\includegraphics[height=3in]{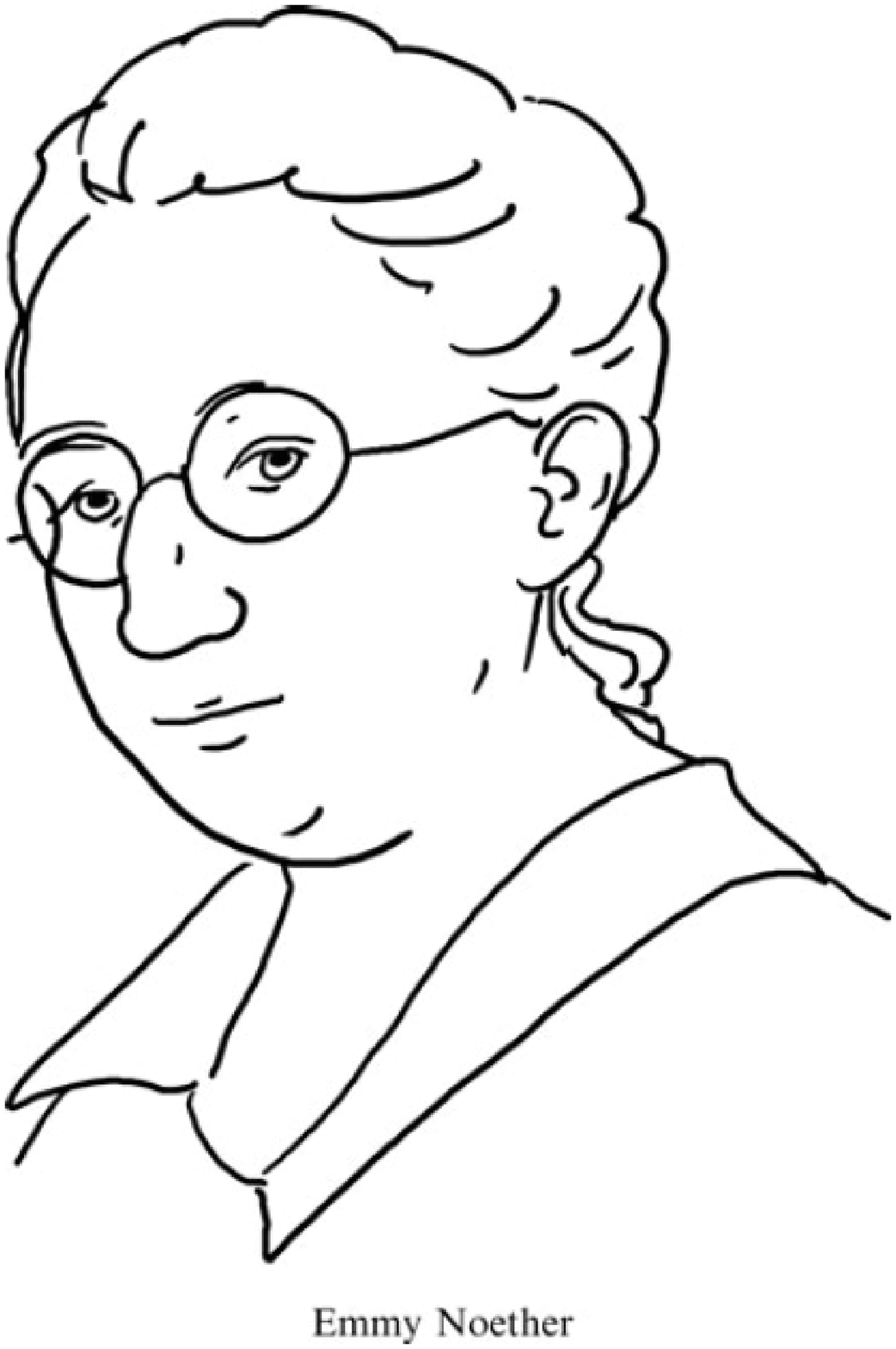}}
\caption{Sketch by G\'eraldine
  d'Alessandris, by permission. {\it Original sketch published in
    ``Vite matematiche'', C. Bartocci et al., Springer, Italy, 2007.}}
\label{fig4}\end{figure}

In other fields of mathematics and in mathematical physics, she is also known
for major contributions. Here, we are concerned with a theorem (now
called ``Noether's theorem''), which she proved in 1915 and published in 1918, that
relates continuous symmetry groups and invariants for mechanical systems
possessing a Lagrangian variational formulation. With the development of
quantum mechanics and field theory in the 20th century, this theorem was to
acquire a central role and is covered in most textbooks on analytical
mechanics or on field theory.\footnote{Noether, 1918.}

It has been known since Lagrange's ``M\'echanique Analitique'' of 1788
that the motion of an incompressible
three-dimensional fluid possesses a variational formulation. In modern 
language, if the fluid occupies the whole space $\R^3$ and the Lagrangian map
is specifed at time $t=0$ (the identity) and at some time 
$T>0$, the Lagrangian map $\bx(\ba,t)$ at intermediate times is an extremum
of the action integral
\UB S = \int_0^Tdt\int_{\R^3} d^3\ba\, L(\dot\bx),\quad {\rm where}\,\,\,L(\dot\bx) \equiv
\frac{1}{2} |\dot\bx|^2, \UE{euleraction}
with the constraint  of incompressibility that 
\UB\det\,(\nL\bx(\ba,t))=1 \quad \text{for all}\,\,\, 0\le t\le T,\UE{sdiffrep}
where we recall that $\nL$ denotes the Lagrangian gradient.
Indeed if, following Lagrange, we introduce infinitesimal variations
$\delta\bx(\ba,t)$, vanishing at $t=0$ and at $t=T$, we find from
\rf{euleraction} that the vanishing of the variation of the action
requires (after an integration by parts over time) that

\UB\delta S =- \int_0^Tdt\int_{\R^3} d^3\ba\, \ddot\bx(\ba,t) \cdot \delta  
\bx(\ba,t) = 0, \UE{lagrangevariation}
for all variations $\delta \bx$ consistent with incompressibility. This 
constraint is more easily written in Eulerian coordinates by defining
\UB \delta ^{\rm E}\bx(\bx,t) \equiv \delta \bx(\ba(\bx,t),t),
\UE{euleriandeltax}
where $\ba(\bx,t)$ is the inverse of the Lagrangian map. The incompressibility
constraint for infinitesimal variations is then simply 
$\nabla \cdot \delta ^{\rm E}\bx(\bx,t) =0$. (In particular, by taking
$\delta \bx = \dot\bx dt$, one has $\nabla\cdot \bv =0$, where $\bv =\dot \bx$.) 
It thus follows that the variation of the action 
\rf{lagrangevariation} must vanish for all $\delta ^{\rm E}\bx$ of zero
Eulerian divergence. Hence the  acceleration $\ddot{\bf x}$ must be the
Eulerian gradient of a suitable function (actually the negative of
the pressure):
\UB \ddot \bx = -\nabla p,\UE{eulermixed}
which is equivalent to the Euler equation \rf{EulerMom} when there is no
external force.

It has been known for a long time that the obvious invariances of
the Lagrangian $(1/2)|\dot\bx|^2$, such as invariance under time and space 
translations and under rotations, are connected, by Noether's theorem, 
to standard mechanical invariants, viz. the conservation of
kinetic energy, of momentum and of angular momentum.  In 1967 William A. Newcomb, a theoretical physicist well-known for the ``Newcomb paradox'',
noticed a new continuous invariance group he called ``exchange invariance''
(which is now mostly called ``relabeling symmetry''). From this he inferred,
by Noether's theorem, new invariants which have been  identified with
the Cauchy invariants many years later 
(see Section~\ref{ss:rebirth}).\footnote{Newcomb, 1967.}
  
Newcomb observed that the action is preserved if  we  change the original Lagrangian coordinates $\ba$
to new Lagrangian coordinates $\ba'$, provided the map from the
 $\ba$\nobreak\hspace{.08em plus .08333em}s to the $\ba'$\nobreak\hspace{.08em plus .08333em}s conserves volumes.  An infinitesimal version, needed
to apply Noethers' theorem, is
\UB \ba \to \ba +\delta \ba(\ba),\qquad\nL\cdot \delta \ba(\ba) =0.\UE{miniexchange} The resulting change in the
Lagrangian map at time $t$ is then (in components with summation over repeated
indices)  $\delta x_i =
\nabla_j^{\rm L} x_i \,\delta a_j$ and the change in the action is
\UB\delta S = \int_0^Tdt\int_{\R^3} d^3\ba\,\, \dot x_i \,\partial_t^{\rm L} 
\left(\nabla_j^{\rm L} x_i\right)\delta a_j, \UE{minichangeaction}
where  $\partial_t^{\rm L}$ denotes the Lagrangian time derivative, when
a dot would be too cumbersome. Setting this variation equal to zero for
all pertubations $\delta \ba$ of vanishing divergence, we find that
$\dot x_i \,\partial_t^{\rm L}  \left(\nabla_j^{\rm L} x_i\right)$  should be
a Lagrangian gradient. Thus its Lagrangian curl should vanish, that is
 \UB \nL \times \left[\partial_t^{\rm L}\left(\dot x_i\nabla_j^{\rm L}
   x_i\right)
 - \ddot x_i  \nabla_j^{\rm L} x_i\right] = 0.\UE{almostthere}
By \rf{eulermixed}, the second term in the square bracket is the Lagrangian
gradient of the pressure, whose Lagrangian curl vanishes. Thus
\rf{almostthere} is equivalent to the Cauchy invariants equations, when they are 
written in Hankel's form, with a curl in front, as in \rf{hankelst}.

Actually, Newcomb was preceded in 1960 by Eckart, a specialist of
quantum mechanics who applied variational methods to fluid
dynamics. Eckart rederived the circulation theorem and obtained the
Cauchy invariants equations (his equations (3.9), (4.4) and (4.5))
without any explicit use of Noether's theorem, but in another 1963
paper Eckart pointed out that ``The general theorems just mentioned
are also consequences of this invariance, a fact that does not seem to
have been noted before.'' (By ``this invariance'' he understands the
unimodular group of volume-conserving transformations of the
Lagrangian coordinates.)  In 1963 Calkin did apply Noether's theorem
to hydrodynamics and magnetohydrodynamics and recovered various known
invariants, but apparently not the Cauchy invariants.  The subject of
the relabeling symmetry was reviewed in Salmon's Annual Review of
Fluid Mechanics paper, which also contains further references, none of
which mentions the Cauchy invariants.\footnote{Eckart, 1960; 1963:
  1038. Calkin, 1963.  Salmon, 1988: \S\,4.}
 
Finally, we should mention that a special case of Cauchy's reformulation
of the equations in Lagrangian coordinates was rediscovered in the
fifties without any use of Noether's theorem. In a classical book on water waves, Stoker pointed out that
certain hydrodynamical problems involving a free surface are better handled
using Lagrangian than Eulerian coordinates. He then described the PhD
work of his student Pohle in the early fifties at the Courant
Institute of Mathematical Sciences in New York. Pohle established
Cauchy's invariant equation \rf{cauchy15} (without citing Cauchy) 
in the special case of two
dimensions when there is a single invariant, while pointing out that
``similar results hold for the three dimensional case.'' He then assumed
analyticity in time of the Lagrangian map and obtained recurrence
relations among the corresponding Taylor coefficients. Here,
two-dimensionality allowed him to use a complex-variable method to
obtain special solutions of relevance to the breakup of a dam. Stoker
stated that the assumed time-analyticity could probably be established
``at least for a finite time'' and pointed out that ``the convergence
of developments  of this kind in some simpler problems in
hydrodynamics has been proved by Lichtenstein.'' This happened indeed, but only
recently. Cauchy's 1815/1827 paper was cited many times in Stoker's
book, in connection with waves, and not in the section discussing
the  Lagrangian formulation. In
the mid-eighties Stoker's 1957 work on the two-dimensional Lagrangian
equations was cited by Abrashkin and Yakubovich, among the persons
who in the late nineties would be involved in correctly identifying
Cauchy's role in the discovery of the invariants (see,
Section~\ref{ss:rebirth}).\footnote{Stoker, 1957: \S\,12.1.  Pohle,
  1951.  Abrashkin \& Yakubovich, 1985.}

\subsection{The renaissance of Cauchy's invariants}
\label{ss:rebirth}

It was only when the 20th century was nearing its end that Cauchy's name was
again associated to his invariants, thanks to Russian scientists. Of
course Russia has had for a long time a very strong tradition in fluid
mechanics. In 1996 Abrashkin, Zen'kovich and Yakubovich, from 
the Institute of Applied Physics in Nizhny Novgorod, wrote a paper
about a new matrix reformulation of the 3D Euler equations in Lagrangian
coordinates. Their eqs.~(4) are Cauchy's invariants equations 
written as three scalar equations, just as in Cauchy's 1815/1827 paper. The time-independent right hand sides are qualified
by them as ``integrals of motion'' and attributed to Cauchy (only by referring
to Lamb). On the next page, they refer to them as ``Cauchy
invariants.'' In a 1997 review paper on nonlinear
waves with significant emphasis on plasma physics, Zakharov and Kuznetsov from
the Landau Institute in Moscow discussed the relabeling
symmetry and the corresponding conservation law. Their eq.~(7.11) gives the
Cauchy invariants equations in vector notation. They do point out, including in their Abstract,  that these are ``the Cauchy invariants.''
Again, they cite  Cauchy's work  with an indirect reference via Lamb.\footnote{Abrashkin, Zen'kovich and Yakubovich, 1996. Lamb, 1932. Zakharov \& Kuznetsov, 1997.}

The name ``Cauchy invariant'' is of course fully appropriate, given
the modern meaning of ``invariant'', a local or global quantity conserved in the
course of time.  This name was already used
in Russia several years prior to the 1996--1997 publications. Evsey Yakubovich (2013 private communication, transmitted
through Evgenii Kuznetsov) confirmed that it was used in internal 
discussions at the  Institute of Applied Physics  in the early nineties. Two
further  publications containing the
name ``Cauchy invariants'' were  published by the Nizhnii
Novgorod group.\footnote{Yakubovich \& Zenkovich, 2001. Abrashkin \& Yakubovich, 2006.}

Cauchy's role in introducing the invariants having thus finally been
recognized, several works discussing the Cauchy invariants
have appeared since the year 2000.\footnote{Friedlander \&
  Lipton-Lifschitz, 2003. Bennett, 2006. Kuznetsov, 2006. Eyink, 2013. 
Ohkitani, 2014. Frisch \& Zheligovsky, 2014. Zheligovsky \& Frisch, 2014.}

\section{Concluding remarks}
\label{s:conclusion}

What have we learned from exploring this more than two-century-long 
(1788--2014) history of the Lagrangian formulation for the 
incompressible Euler equations?  Actually, we have a situation (here
in hydrodynamics), for which a discovery made
two centuries back and hardly ever used since, has emerged as particularly
relevant for modern research developments.  We have in mind Pohle's early-fifties work on the
breakup of dams, the invariants obtained from Noether's theorem in the
sixties and 
the recent proof of time-analyticity of Lagrangian fluid
particle trajectories. In all these instances, Cauchy's invariants
equations do play a key role, but were actually rediscovered in
different ways, for example by a  transposition to
incompressible hydrodynamics of a Lagrangian perturbation method that
cosmologists have developed since the nineties.\footnote{On the relation of Cauchy's
invariants equations and cosmology, see Zheligovsky \& Frisch, 2014:
\S\,3.}

Probably, Cauchy made
his 1815 discovery too early to be adequately appreciated, because the study
of invariants
and conservation laws would not emerge as an important paradigm for
another century. Only a 
corollary of \rf{cauchy15}, namely Cauchy's vorticity formula
\rf{cauchy17} was attracting attention in the early 19th century, because of 
Lagrange's theorem. The first serious opportunity to
understand some of the importance of Cauchy's invariants equations came
in Germany around 1860 after Dirichlet stated that Lagrangian coordinates
are important (a statement not heard much again until late in
the 20th century), when G\"ottingen University, possibly prompted by Riemann,
pushed for studies based on Lagrangian
coordinates to get more insight into Helmholtz's vorticity theorems and,
last but not least, when Hankel found how to make full use of  Cauchy's
invariants equations and in the process came across the circulation
theorem some years before Kelvin. 

Later, in the first decades  of the
20th century, Cauchy's invariants equations faded into oblivion or were
viewed a mere intermediate step in proving  Cauchy's vorticity
formula. Eventually, developments in theoretical mechanics, closely
connected to the rise of quantum mechanics and later of quantum
field theory, led to a rediscovery of the invariants in the late 1960 
through application of Noether's theorem. Another 30 years elapsed, during
which developments in nonlinear physics  were increasingly making use of symmetries and
invariants, until the crucial importance of Cauchy's work could be
appreciated.\\

\noindent ACKNOWLEDGMENTS 

\noindent We are grateful to J\'er\'emie Bec, Bruno Belhoste, Olivier Darrigol, Michael Eckert, 
Evgenii Kuznetsov, Koji Ohkitani, Walter Pauls and Vladislav Zheligovsky for useful 
discussions and to Florence Greffe for helping us with historical material.

\end{document}